\documentclass[12pt]{amsart}
\usepackage{latexsym, amsmath, amscd, amssymb, amsthm}

\newtheorem*{te}{Theorem}
\newtheorem*{lm}{Lemma}

\begin{document}

\noindent

 \title{  On complete system of covariants for the   binary form of  degree 7}

\author{Leonid Bedratyuk} \address{ Khmelnytsky National University, Instytuts'ka st. 11, Khmelnytsky , 29016, Ukraine}
\email {bedratyuk@ief.tup.km.ua}
\begin{abstract}
A minimal system of 147 homogeneous generating elements of the algebra of covariants for the binary form of degree 7 is calculated.
\end{abstract}

\maketitle

\section{introduction}

Let $V_n$ be a vector $\mathbb{C}$-space of the binary forms of degree $d$ considered with natural action of the group  $G=SL(2,\mathbb{C}).$ Let us extend the action of the group  $G$ to the polynomial functions algebra  $\mathbb{C}[V_d \oplus \mathbb{C}^2 ].$ 
Denote  by  ${C_d=\mathbb{C}[V_d \oplus \mathbb{C}^2 ]^{\,G}}$ the corresponding subalgebra of   $G$-invariant functions. 
In the vocabulary of classical invariant theory  the algebra    $C_d$ is called     the  algebra covariants of the  binary form of d-th degree.
Let  $C_d^{+}$ be an ideal of  $C_d$ generated by all homogeneus elements of positive power. Denote by  $\bar C_d$ a set of  homogeneus elements of $C_d^{+}$ such that their images in  $C_d^{+}/(C_d^{+})^2$ form a basis of the vector space. The  set  $\bar C_d$ is caled complete system of covariants of the  d-th degree binary form. Elements of   $\bar C_d$ form a minimal system of homogeneous generating elements of the invariants algebra $C_d.$ Denote by $c_d$ a number of  elements of the set    $\bar C_d.$

The complete systems of  covariants was a topic of major research interest in classical invariant theory of the 19th century.
It is easy to show that $c_1=0, $  $c_2=2,$ $c_3=4.$ A complete system of covariants  in the case   $d=4$ was calculated by Bool, Cayley, Eisenstein, see survey \cite{Dix}.
The complete systems  of invariants and covariants in the cases   $d=5,6$ were calculated by Gordan, see  \cite{Gor}. In particular, $c_4=5,$  $c_5=23,$ $c_6=26.$

Gall's attempt  \cite{Gall-1}  to discover the complete system for the 
 case $d=7$ was unsuccessful. He did offered a system of 151 covariants but 
 the system was not a minimal system, see \cite{DL},\cite{Aut-1}. Also, 
 Sylvester's attempts \cite{SF} to find the cardinality 
 $c_7$ and covariant's degree-order distribution  of the complete system 
 were mistaken, see \cite{DL}.  Therefore, the problem of finding a minimal system of 
 homogeneous generating elements (or even a cardinality of the system) of 
 the  algebra  of covariants for  the  binary form of degree 7 was  
 open.

 The case  $d=8$  was studied by Sylvester and Gall but they have obtained 
 completelly different results. By using Sylvester-Cayley technique, 
 Sylvester in   \cite{SF} got that  ${c_8=69.}$   Gall in    \cite{Gall}, 
 evolving the Gordan's constructibe method, offered   68 covariants as a 
 complete system of covariants for  the case  $d=8.$  In \cite{Aut-8} the autor have calculated the set $\bar C_8.$ 

 For the cases $d=9,10$ Sylvester in \cite{SF} calculated the 
 cardinalities  $c_9$ and  $c_{10}$  but the present author, using a 
 computer, found numerous mistakes in those computations.

 Therefore, the 
 complete systems of covariants for the binary form are so far 
 known only for $d\leq 6$, see \cite{Dix}, \cite{Olver},  and  for $d=8$, see \cite{Aut-8}.

This article's purpose is calculating the set $\bar C_7.$  Using a technique introduced in \cite{Aut-1}  we explicitly calculate a complete system of covariants of the 7th degree binary form. In particular, we prove that $c_7=147.$

All calculation were done with Maple.
\section{Premilinaries. }

 
Before any calculation we try make a simplification of a covariant represenation and their computation. 

Let us indentify the algebra $\mathbb{C}[V_d]$  with the algebra $\mathbb{C}[X_d]:=\mathbb{C}[t,x_1,x_2,\ldots, x_d],$ and the algebra $\mathbb{C}[V_d \oplus \mathbb{C}^2 ]$ identify  with the polynomial algebra $\mathbb{C}[t,x_1,x_2,\ldots, x_n,Y_1,Y_2].$
Let  $${\it \kappa} : C_d \longrightarrow \mathbb{C}[X_d]$$ be the   $\mathbb{C}$-linear map takes each homogeneous covariant of order $k$ to his leading coefficient, i.e. a coefficient of   $Y_1^k.$ Follow by classical  tradition an element of the algebra   ${\it \kappa} (C_d)$  is called  {\it semi-invariants, }  a degree of a homogeneous covariant with respect to the variables set $X_d$ is called {\it degree} of the covariant and its  degree with respect to the variables set $Y_1, Y_2$  is called {\it order.}

In \cite{Aut-1} showed that 
$$
{\it \kappa} (C_d)=k[t,z_2, z_3,\ldots ,z_d][\frac{1}{t}]\cap \mathbb{C}[X_d], 
$$
$$
z_i:= \sum_{k=0}^{i-2} (-1)^k {i \choose k} x_{i-k}  x_1^k t^{i-k-1} +(i-1)(-1)^{i+1} x_1^i, i=2,\ldots,d.
$$
Let us consider the following derivation $D$ of the algebra  $k[t,x_1,z_2, z_3,\ldots ,z_d][\frac{1}{t}]:$ 

$$ 
\begin{array}{l}

D=\displaystyle 7\,{x_{1}}\,{\frac {\partial }{\partial t}} - 
\displaystyle \frac {( - 15\,{x_{1}}\,{z_{3}} + 18\,{z_{2}}^{2}
 - 4\,{z_{4}})}{t}{\frac {\partial }{\partial {z_{3}}}} 
 + {\displaystyle \frac {(20\,{x_{1}}\,{z_{4}} - 24\,{z_{2}}\,{z
_{3}} + 3\,{z_{5}})}{
t}}\,{\frac {\partial }{\partial {z_{4}}}} +  \\
\\
\mbox{} + {\displaystyle \frac {(2\,{z_{6}} + 25\,{x_{1}}\,{z_{5}
} - \displaystyle 30\,{z_{2}}\,{z_{4}})}{t}}\,{\displaystyle \frac {\partial }{\partial {z_{5}}}}
  + {\displaystyle \frac {({z_{7}} + 30\,{x_{1}}\,{z_{6}
} - 36\,{z_{2}}\,{z_{5}})}{t}} \,{\displaystyle \frac {\partial }{\partial {z_{6}}}}+
 \\
\\
\mbox{} + {\displaystyle \frac {7\,(5\,{x_{1}}\,{z_{7}} - 6\,{z_{
2}}\,{z_{6}})}{t}} \,{\displaystyle \frac {\partial }{\partial {z_{7}}}}
 + {\displaystyle \frac {5\,(2\,{x_{1}}\,{z_{2}} + {z_{3}})}{t}} \,
{\displaystyle \frac {\partial }{ \partial {z_{2}}}}.
\end{array}
$$

Suppose  $F=\displaystyle \sum_{i=0}^m \, f_i { m \choose i } Y_1^{m{-}i} Y_2^i$ be a covariant of order $m,$  $\kappa(F)=f_0 \in {\it \kappa} (C_d).$
The classical Robert's theorem \cite{Rob},  states that the covariant $F$ is completelly and   uniquely determined  by its leading coefficient $f_0,$ namely, see  \cite{Aut-1} :
$$
F=\sum_{i=0}^{m} \frac{D^i(f_0)}{i!} Y_1^{m-i}Y_2^i.
$$

On the other hand, every semi-invariant is a leading coefficient of a  covariant, see \cite{Gle}, \cite{Olver}. This 
give us well defined   explicit form of the inverse map  $${\kappa^{-1} :   {\it \kappa} (C_d) \longrightarrow C_d,}$$ namely

$$
\kappa^{-1}(a)=\sum_{i=0}^{{\rm ord}(a)} \frac{D^i(a)}{i!} Y_1^{{\rm ord}(a)-i}Y_2^i,  
$$
here  $ a \in   {\it \kappa} (C_d)$ and  ${\rm ord}(a)$ is an order of the element $a$ with respect to the locally nilpotent (at ${\it \kappa} (C_d)$) derivation  $D,$  i.e. ${\rm ord}(a):=\max \{ s, D^s(a) \neq 0 \}.$  For example, since ${\rm ord}(t)=d,$ we have 
$$
\kappa^{-1}(t)=\sum_{i=0}^{{\rm ord}(t)} \frac{D_2^i(t)}{i!} Y_1^{{\rm ord}(t)-i}Y_2^i =t Y_1^d+\sum_{i=1}^{d} { d \choose i } x_i Y_1^{d-i}Y_2^i.
$$
As we see  the $\kappa^{-1}(t)$ is just the  basic binary form. From polynomial functions point of view the covariant $\kappa^{-1}(t)$ is the evaluation  map.

Thus,  the  problem of finding of complete system of the  algebra   $\overline C_d$ is equivalent to the  problem  of finding of complete system of semi-covariants's algebra   $\kappa(C_d).$  It is well known classical results, see, for example, \cite{Cay}.

To calculate the semi-invariants we need have an analogue of the transvectants. 
Let  $f, g$ are two covariants of the degrees  $m$ and $ k$ respectively. The semi-invariant of the form  

$$
[f,g]^r:=q_r(f,g) \cdot \sum_{i=0}^r (-1)^i { r \choose i } \frac{D^i(f)}{[m]_i} {\Bigl | _{x_1=0,\ldots,x_r=0}}  \frac{D^{r-i}(g)}{[k]_{r-i}}{\Bigl | _{x_1=0,\ldots,x_r=0}}, \mbox{ } 0 \leq r \leq \min(m,k),
$$
is called a $r$-th semintransvectant of the semi-invariants $f$ and $g.$
Here  $q_r(f,g)$ is a normalization rational factor, see \cite{Aut-1}. Note that up to multiplication on a rational number the semi-transvectant  $[f,g]^r$ is equal to  $\kappa(\kappa^{-1}(f), \kappa^{-1}(g))^r, $ where  $(\cdot , \cdot )^r$   is usual transvectant,  see  \cite{Hilb}, \cite{Olver}.
The following statements is direct consequences of corresponding  transvectant properties,  see \cite{Gle}:
\begin{lm} Let  $f, g$ be two semi-invariants.  Then the folloving conditions hold
\begin{enumerate}
\item[({\it i})] the semitransvectant $[t, f\,g]^i$ is reducible for  $ 0 \leq i \leq \min(d, \max( {\rm ord}(f),{\rm ord}(g));$
\item[({\it ii})] if  $ {\rm ord}(f)=0,$ then   $[t, f\,g]^i=f [t,g]^i;$
\item[({\it iii})] ${\rm ord}([f,g]^i)={\rm ord}(f)+{\rm ord}(g)-2\, i;$ 
\item[({\it iv})] ${\rm ord}(z_2^{i_1}z_3^{i_3} \cdots z_d^{i_d})=d\,(i_2+i_3+\cdots +i_d)-2\,(2\, i_2+3\, i_3+\cdots +d \, i_d).$
\end{enumerate}
\end{lm}

The representation of semi-invariants as elements of the algebra  
  $\displaystyle \mathbb{C}[t,z_2,z_3,\ldots z_d][\frac{1}{t}]$  is more compact than their standard representation as elements of the algebra  $\mathbb{C}[X_d].$ Very rough empire estimate is  that a semi-invariant $\kappa(F)$ has terms number  in ${\displaystyle ([\deg (F)/d]+2)\, ({\rm ord}(F)}+1)$ times less  than the corresponding covariant $F.$ Moreover, from computing  point of view, the semitransvectant formula is more effective than the transvectant  formula. These two favorable circumstances coupled with great Maple powers allow us  to compute  a complete system of covariants of the  7-th degree binary form.

\section{Irreducible covariants up to 13th degree}


Let  $\overline C_{7,\,i}:= (C_7)_i$ be a subset  of  $\overline C_{7}$ whose elements has degree  $i$ 
 and  $\overline C_{7,\,i}:=\overline C_7 \cap C_{7,\,i}.$
By using an analogue of the $\Omega$-process and the lemma, in \cite{Aut-1} was found all irreducible covariants up to 13th degree. Let us present the lists of the generating elements.

The unique semi-invariant of the degree one obviously is  $t, {\rm ord}(t)=7.$

The set  $\overline C_{7,\,2}, $ consists of the 3  irreducible semi-invariants 
$$
\begin{array}{ll}
dv_1:=[t,t]^4={\displaystyle \frac {3\,{z_{2}}^{2} + {z_{4}}}{t^{2}}}={x_{4}}\,t - 4\,{x_{1}}\,{x_{3}} + 3\,{x_{2}}^{2} , & {\rm ord}(dv_1){=}6,  \\ dv_2:=[t,t]^6={\displaystyle \frac {{z_{6}} + 15\,{z_{2}}\,{z_{4}} - 10\,{z_{3}, }^{2}}{t^{4}}}={x_{6}}\,t - 6\,{x_{1}}\,{x_{5}} + 15\,{x_{2}}\,{x_{4}} - 10\,{x
_{3}}^{2}, & {\rm ord}(dv_2){=}2, \\
dv_3:=[t,t]^2=z_2={x_{2}}\,t - {x_{1}}^{2}, &  {\rm ord}(dv_2){=}10.
\end{array}
$$

The set  $\overline C_{7,\,3}, $ consists of the 6 irreducible semi-invariants
$$
\begin{array}{llll}
tr_1=[t,dv_1]^4, & {\rm ord}(tr_1)=5, & tr_2=[t,dv_3],  & {\rm ord}(tr_2)=15,\\

tr_3=[t,dv_3]^3, &  {\rm ord}(tr_3)=11, & tr_4=[t,dv_3]^4, & {\rm ord}(tr_4)=9,\\

tr_5=[t,dv_3]^5,  & {\rm ord}(tr_5)=7, & tr_6=[t,dv_3]^7, & {\rm ord}(tr_6)=3.

\end{array}
$$

The set $\overline C_{7,\,4}$ consists of the following 8 ireducible semi-invariants: 
$$
\begin{array}{llll}
ch_1=[t,tr_5]^7, & {\rm ord}(ch_1)=0, & ch_2=[t,tr_3]^7,  & {\rm ord}(ch_2)=4,\\

ch_3=[t,tr_3]^2, &  {\rm ord}(ch_3)=14, & ch_4=[t,tr_3]^4,  & {\rm ord}(ch_4)=10,\\

ch_5=[t,tr_3]^5, &  {\rm ord}(tr_5)=8, & ch_6=[t,tr_1]^2, & {\rm ord}(tr_6)=8,\\

ch_7=[t,tr_1]^3, & {\rm ord}(ch_7)=6, & ch_8=[t,tr_1]^4, & {\rm ord}(ch_6)=4.
\end{array}
$$

The set $\overline C_{7,\,5}$ consists of the following 10 ireducible semi-invariants:
$$
\begin{array}{llll}
pt_1=[t,ch_6]^5, &{\rm ord}(pt_1)=5, & pt_2=[t,ch_6]^6, & {\rm ord}(pt_2)=3,\\

pt_3=[t,ch_7]^2, & {\rm ord}(pt_3)=9, & pt_4=[t,ch_7]^3,  & {\rm ord}(pt_4)=7,\\

pt_5=[t,ch_7]^5,  & {\rm ord}(pt_5)=3, & pt_6=[t,ch_6]^3,  & {\rm ord}(pt_6)=9,\\

pt_7=[t,ch_4]^2,  & {\rm ord}(pt_7)=13, & pt_8=[t,ch_4]^5, & {\rm ord}(pt_6)=7, \\

pt_9=[t,dv^2_1]^7,  & {\rm ord}(pt_9)=5, & pt_{10}=[t,dv_1 dv_2]^7,  & {\rm ord}(pt_{10})=1. \\
\end{array}
$$

The set $\overline C_{7,\,6}$ consists of the following 10 ireducible semi-invariants:

$$
\begin{array}{llll}
sh_1=[t,pt_5]^5, & {\rm ord}(sh_1)=6, & sh_2=[t,pt_7]^6, &  {\rm ord}(sh_2)=8,\\

sh_3=[t,pt_4]^5, &  {\rm ord}(sh_3)=4, & sh_4=[t,pt_4]^6, & {\rm ord}(sh_4)=2,\\

sh_5=[t,pt_3]^2, &  {\rm ord}(sh_5)=12, & sh_6=[t,pt_3]^4,  & {\rm ord}(sh_6)=8,\\

sh_7=[t,pt_4]^4,  & {\rm ord}(sh_7)=6, & sh_8=[t,tr_1 dv_1]^7, & {\rm ord}(sh_6)=4, \\

sh_9=[t,tr_1 dv_2]^6,  & {\rm ord}(sh_9)=2, & sh_{10}=[t, tr_6 dv_1]^7, &  {\rm ord}(sh_{10})=2. \\
\end{array}
$$
       
 The set $\overline C_{7,\,7}$ consists of the following 12 ireducible semi-invariants: 

$$
\begin{array}{llll}
si_1=[t,sh_5]^4, & {\rm ord}(si_1)=11, & si_2=[t,sh_7]^4, & {\rm ord}(si_2)=5,\\

si_3=[t,tr_1^2]^7,  & {\rm ord}(si_3)=3, & si_4=[t,sh_1]^3, & {\rm ord}(si_4)=7,\\

si_5=[t,ch_7 dv_1]^7,  & {\rm ord}(si_5)=5, & si_6=[t,ch_7 dv_2]^7,  & {\rm ord}(si_6)=1,\\

si_7=[t,tr_6^2]^4, & {\rm ord}(si_7)=5, & si_8=[t,tr_6^2]^6,  & {\rm ord}(si_6)=1, \\

si_9=[t,tr_6\, tr_1]^6,  & {\rm ord}(si_9)=3, & si_{10}=[t, tr_6\, tr_1]^7, & {\rm ord}(si_{10})=1, \\

si_{11}=[t,tr_1^2]^6,  & {\rm ord}(si_{11})=5, & si_{12}=[t,sh_{10}), & {\rm ord}(si_{12})=7.
\end{array}
$$      

   The set $\overline C_{7,\,8}$ consists of the following 13 ireducible semi-invariants:     
$$
\begin{array}{llll}
vi_1=[t,si_7]^3, & {\rm ord}(vi_1)=6, & vi_2=[t,si_7]^4,  & {\rm ord}(vi_2)=4,\\

vi_3=[t,ch_8\,tr_6]^7,  & {\rm ord}(vi_3)=0, & vi_4=[t,ch_8\,tr_1]^6,  & {\rm ord}(vi_4)=4,\\

vi_5=[t,ch_8 tr_1]^7,  & {\rm ord}(vi_5)=2, & vi_6=[t,ch_7 tr_6]^7,  & {\rm ord}(vi_6)=2,\\

vi_7=[t,ch_7 \,tr_1]^7,  & {\rm ord}(vi_7)=4, & vi_8=[t,ch_8\,tr_6]^6, & {\rm ord}(vi_6)=2, \\

vi_9=[t,tr_6\, dv_2^2]^7, & {\rm ord}(vi_9)=0, & vi_{10}=[t,si_4]^2, & {\rm ord}(vi_{10})=10, \\

vi_{11}=[t,si_12]^4, & {\rm ord}(vi_{11})=6, & vi_{12}=[t,si_{11}]^3),  & {\rm ord}(vi_{12})=6,\\
& vi_{13}=[t,pt_9\,dv_2]^7, & {\rm ord}(vi_{13})=0.&
\end{array}
$$

The set $\overline C_{7,\,9}$ consists of the following 11 ireducible semi-invariants:
$$
\begin{array}{llll}
de_1=[t,sh_3\,dv_1]^7, & {\rm ord}(de_1)=3, & de_2=[t,ch_7\,ch_8]^7, & {\rm ord}(de_2)=3,\\

de_3=[t,pt_5\,tr_6]^5, & {\rm ord}(de_3)=3, & de_4=[t,pt_5 \,tr_1]^6,  & {\rm ord}(de_4)=3,\\

de_5=[t,pt_5\, tr_1]^7,  & {\rm ord}(de_5)=1, & de_6=[t,sh_9\,dv_1]^7,  & {\rm ord}(de_6)=1,\\

de_7=[t,sh_{10}\,dv_1]^7,  & {\rm ord}(de_7)=1, & de_8=[t,sh_{10} \,dv_2]^3, & {\rm ord}(de_6)=5, \\

de_9=[t,vi_5]^2,  & {\rm ord}(de_9)=5, & de_{10}=[t,vi_2]^4, & {\rm ord}(de_{10})=3, \\

& de_{11}=[t,vi_{11}]^2, & {\rm ord}(de_{11})=9.&
\end{array}
$$

The set $\overline C_{7,\,10}$ consists of the following 9 ireducible semi-invariants:
$$
\begin{array}{llll}
des_1=[t,sh_9\,tr_1]^6, & {\rm ord}(des_1)=2, & des_2=[t,sh_4\,tr_6]^4, & {\rm ord}(des_2)=4,\\

des_3=[t,sh_4\,tr_1]^6,  & {\rm ord}(des_3)=2, & des_4=[t,sh_1\,tr_1]^7, & {\rm ord}(des_4)=4,\\

des_5=[t,sh_3\,tr_6]^5, & {\rm ord}(des_5)=4, & des_6=[t,de_9]^2, & {\rm ord}(des_6)=8,\\

des_7=[t,tr_6^3]^7, & {\rm ord}(des_7)=2, & des_8=[t,sh_{10} \,tr_1]^6, & {\rm ord}(des_6)=2, \\

& des_9=[t,pt_1\,ch_7]^7,  & {\rm ord}(des_9)=4. &
\end{array}
$$

The set $\overline C_{7,\,11}$ consists of the following 9 ireducible semi-invariants:
$$
\begin{array}{llll}
odn_1=[t,vi_2\,dv_1]^7, &{\rm ord}(odn_1)=3, & odn_2=[t,vi_2,dv_2]^6,  & {\rm ord}(odn_2)=1,\\

odn_3=[t,vi_4\, dv_2]^6, & {\rm ord}(odn_3)=1, & odn_4=[t,vi_5\,dv_1]^7, & {\rm ord}(odn_4)=1,\\

odn_5=[t,vi_6\, dv_1]^7, & {\rm ord}(odn_5)=1, & odn_6=[t,vi_2\,dv_2]^5,  & {\rm ord}(odn_6)=3,\\

odn_7=[t,des_6]^4, & {\rm ord}(odn_7)=7, & odn_8=[t,des_6]^6,  &{\rm ord}(odn_6)=3, \\

& odn_9=[t,vi_1\,dv_2]^7, & {\rm ord}(odn_9)=1.&
\end{array}
$$

The set $\overline C_{7,\,12}$ consists of the following 13 ireducible semi-invariants:
$$
\begin{array}{llll}
dvan_1=[t,sh_1\,pt_2]^7, & {\rm ord}(dvan_1)=2, & dvan_2=[t,sh_1\,pt_5]^7, & {\rm ord}(dvan_2)=2,\\

dvan_3=[sh_9,sh_{10}]^2, & {\rm ord}(dvan_3)=0, & dvan_4=[t,odn_7]^6, & {\rm ord}(dvan_4)=2,\\

dvan_5=[t, de_8\,dv_2]^6, & {\rm ord}(dvan_5)=2, & dvan_6=[sh_{10}\, ,sh_{10}]^2, & {\rm ord}(dvan_6)=0,\\

dvan_7=[t, de_9\,dv_2]^6, & {\rm ord}(dvan_7)=2, & dvan_8=[t, de_{10}\,dv_1]^7,  & {\rm ord}(dvan_6)=2, \\

dvan_9=[t,odn_7]^4, & {\rm ord}(dvan_9)=6, & dvan_{10}=[sh_1\, , sh_1]^2, & {\rm ord}(dvan_{10})=0, \\

dvan_{11}=[sh_4\, , sh_4]^2, & {\rm ord}(dvan_{11})=0, & dvan_{12}=[sh_4, sh_9]^2), & {\rm ord}(dvan_{12})=0,\\

&  dvan_{13}=[sh_4\, , sh_2]^2, & {\rm ord}(dvan_{13})=0. &
\end{array}
$$   
  
The set $\overline C_{7,\,13}$ consists of the following 9 ireducible semi-invariants:
$$
\begin{array}{llll}
tryn_1=[t,dvan_9]^6, & {\rm ord}(tryn_1)=1, & tryn_2=[t,vi_1\,ch_7]^7, &  {\rm ord}(tryn_2)=5,\\

tryn_3=[t,vi_2\,ch_8]^7, & {\rm ord}(tryn_3)=1, & tryn_4=[t,vi_2\,ch_2]^7, & {\rm ord}(tryn_4)=1,\\

tryn_5=[t,vi_1\,ch_8]^7, & {\rm ord}(tryn_5)=3, & tryn_6=[t,vi_5\,ch_2]^6, & {\rm ord}(tryn_6)=1,\\

tryn_7=[t,vi_8\,ch_8]^6, & {\rm ord}(tryn_7)=1, & tryn_8=[t,vi_8\,ch_7]^7,  & {\rm ord}(tryn_6)=1, \\

& tryn_9=[t,vi_4\,ch_8]^7, & {\rm ord}(tryn_9)=1. &

\end{array}
$$  

The cardinalities of the sets  $\overline C_{7,\,i}$ for $i=1,\ldots,13$  and the covariant's degree-order distributions so far coincide  completelly with  Gall's  and Sylvester's    results, see \cite{Gall-1}, \cite{SF} .

  \section{Computation of covariants of  degree 14--30}

Let  $(C_{+}^2)_i$ be a subset  of  $(C_7^{+})^2$ whose elements  has degree  $i.$ Denote by $C_{i,j}$  a subspase of $C_{7,i}$ generated by elements of degree $i$ and order $j.$ 
Denote by  $\delta_i$ a  number of  irreducible invariants of degree $i$ and denote by  $\delta_{i,j}$ a number of  irreducible invariants of degree $i$ and order $j.$ It is evident that $\delta_i=\sum_{j} \delta_{i,j}$ and $c_7=\sum_i \delta_i.$
To compute the number $\delta_{i,j}$ we use the formula   $\delta_{i,j} =\dim C_{i,j}-\dim(C^2_{+})_{i,j}.$ A dimension of the vector space  $C_{i,j}$ is calculated by Cayley-Sylvester formula, see, for example  \cite{Hilb}, \cite{SP} and  the dimension is equal to a coefficient of 
$\displaystyle T^{\frac{d*i-j}{2}}$ in the expansion of the series 
$$
\frac{(1-T^{d+1})(1-T^{d+2})\ldots (1-T^{d+i})}{(1-T^2)\ldots(1-T^i)}.
$$ 

A dimension of the vector subspace  $(C^2_{+})_{i,j}$ of $(C^2_{+})_{i}$ whose elements has order $j$  is calculated by the formula $\dim(C^2_{+})_{i,j}=\sigma_{i,j} -\dim S_{i,j}.$ Here  $\sigma_{i,j}$ is a number of monomial of the vector space  $(C^2_{+})_{i,j}$  and  $S_{i,j}$ is a vector space of  $(C^2_{+})_{i,j}$ generated by syzygies. A dimension of the vector spase $S_{i,j}$ is found by direct Maple computations. A calculation  of the numbers $\sigma_{i,j}$ is a simple combinatorial problem. An explicit way of  the  irreducible semi-invariants is found by linear algebra techinque, see details in \cite{Aut-1}, \cite{Aut-8}.

The set   $\overline C_{14,0}$  we take from  \cite{Aut-1}: $$chot_{1}:=[si_8,si_{10}], chot_2:= [si_6,si_{10}], chot_3:=[si_6,si_8], chot_4:=[si_3,si_9]^3.$$
 For $j=2$ we have  $\dim C_{14,2}=30,$ $\sigma_{14,2}=36,$ $\dim S_{14,2}=6.$ Thus  $\delta_{14,2}=30-(36-6)=0.$
 For $j=4$ we have  $\dim C_{14,4}=37,$ $\sigma_{14,4}=60,$ $\dim S_{14,4}=25.$ Thus  $\delta_{14,4}=37-(60-25)=2.$  After calculation we get  the following irreducible semi-invariants of   $\overline C_{14,2}:$ 
$$chot_5:=[sh_{10},vi_2],  chot_6:=[sh_9,vi_2].$$
The results of Sylvester's and Gall's  for $\delta_{14,4}$ were 0 and 2 respectively.
Taking into account the results of \cite{Gall-1} we have  $\delta_{14,j}=0$ for $j>4$. Therefore $\delta_{14}=6$  and the set $C_{7,14}$  consists of the following  6 semi-invariants -- $chot_1, \ldots , chot_6.$

There exists  only  semi-invariants of orders 1 and 3 in the set  $\overline C_{15},$ see \cite{Gall-1}. We have  ${\dim C_{15,1}=20,}$ $\sigma_{15,1}=17,$ $\dim S_{15,1}=0.$ Consequently $\delta_{15,1}=3.$ In the same way we get   $\dim C_{15,3}=42,$ $\sigma_{15,3}=61,$ $\dim S_{15,3}=20.$ Thus $\delta_{15,3}=1$ and  $\delta_{15}=4.$
The results of Sylvester's and Gall's $\delta_{15,3}$ were 0 and 1 respectively.
After calculation we obtain the 4 irreducible semi-invariants:
$$ptn_{1}:= [de_{10},sh_9]^2, ptn_{2}=[de_{10},sh_4]^2, ptn_{3}=[de_{3},sh_9]^2, ptn_{4}=[de_{10},sh_{10}].$$

There exists  only  semi-invariants of orders 0 and 2 in the set  $\overline C_{16}.$ The two irreducible invariants were found in  \cite{Aut-1} -
$$shis_{1}:=[vi_2,vi_4]^4, shis_{2}:=[vi_4,vi_7]^4. $$
For $j=2$ we have $\dim C_{16,2}=33,$ $\sigma_{16,2}=39,$ $\dim S_{16,2}=9.$ Thus $\delta_{16,2}=3 $  and we have $\delta_{16}=5.$ 
The results of Sylvester's and Gall's  for $\delta_{16,2}$ were 0 and 3 respectively.
After calculation we obtain the rest 3 semi-transvectants of the order 2:
$$
shis_{3}= [vi_{5},vi_{2}]^2, shis_{4}= [vi_{8},vi_{2}]^2, shis_{5}= [des_{7},sh_{10}].
$$

There exists only semi-invariants of  order  1 in $\overline C_{17}$. We have  $\dim C_{17,1}=31,$ $\sigma_{17,1}=29,$ $\dim S_{17,1}=0.$  Therefore  $\delta_{17}=\delta_{17,1}=2. $ 
Sylvester and Gall have got the same value for $\delta_{17,1}.$
After calculation we obtain the 2 irreducible semi-invariants:
$$
simn_{1}:=[de_{3},vi_{5}]^ 2, simn_{2}:=[si_{8},des_{7}]. 
$$

There exists  only  semi-invariants of orders 0 and 2 in the set  $\overline C_{18}.$ The nine  irreducible invariants were found in  \cite{Aut-1} -

$$
\begin{array}{l}
vis_{1}:=[de_4,de_3]^ 3, vis_{2}:=[de_4,de_{10}]^ 3, vis_{3}:=[de_5,de_6], vis_{4}:=[de_1,de_{10}]^ 3, vis_{5}:=[de_2,de_3]^ 3, \\
\\
vis_{6}:=[de_2,de_{10}]^ 3, vis_{7}:=[de_3,de_{10}]^ 3, vis_{8}:=[de_6,de_{7}],  vis_{9}:=[de_8,de_{9}]^5.
\end{array}
 $$
For  $\overline C_{18,2}$ we have   $\dim C_{18,2}=63,$ $\sigma_{18,2}=105,$ $\dim S_{18,2}=42.$ Thus  $\delta_{18,2}=0 $ i $\delta_{18}=9.$
The results of Sylvester's and Gall's  for $\delta_{18,2}$ were 0 and 1 respectively.

There exists only semi-invariants of  order  1 in $\overline C_{19}$. 
We have  $\dim C_{19,1}=46,$ $\sigma_{19,1}=57,$ $\dim S_{19,1}=12.$ Thus  $\delta_{19}=\delta_{19,1}=1.$ 
The results of Sylvester's and Gall's  for $\delta_{19,1}$ were 0 and 2 respectively.
After calculation we obtain the  irreducible semi-invariant: 
$$
devn:=[de_7,des_7].
$$ 

There exists only one invariant  in $\overline C_{20}$ which is calculated in 
  \cite{Aut-1} -- ${dvad_{20}:=[des_7,des_7]^2.}$
The results of Sylvester's, Gall's   and Diximier's  for $\delta_{20}$ were 0, 2 and 1 respectively.

There exists only two  invariants  in $\overline C_{22}$ which are  calculated in 
  \cite{Aut-1}:
 $$dvdv_{1}:=[odn_6, odn_1]^3, dvdv_{2}:=[odn_8, odn_1]^3.$$
The results of Sylvester's, Gall's   and Diximier's  for $\delta_{22}$ were 1, 3 and 2 respectively.

There exists only semi-invariants of  order  1 in $\overline C_{23}$.
We have   $\dim C_{23,1}=85,$ $\sigma_{23,1}=142,$ $\dim S_{23,1}=58.$ Thus  $\delta_{23}=\delta_{23,1}=1.$
The results of Sylvester's and   Gall's  for $\delta_{23}$ were 0  and 1 respectively.
After calculation we get the semi-invariant -- 
$dvtr:=[tryn_4,des_7].$

In  \cite{Gall-1}, \cite{DL}, \cite{Aut-1} proved  that $\delta_{24}=0.$

There exists only semi-invariants of  order  1 in $\overline C_{25}$.
We have   $\dim C_{25,1}=114,$ $\sigma_{25,1}=228,$ $\dim S_{25,1}=114.$ Òîìó $\delta_{25}=\delta_{25,1}=0.$
The results of Sylvester's and   Gall's  for  $\delta_{25}$ were 0  and 1 respectively.

There exists only one covariant  of  degree 26 -- $dvsh=[tryn_4,tryn_3]$, no  covariants of degrees $27,28,29$  and only one covariant of degre 30  --  $ trd:=(h,\alpha) $ (in Gall's notations).
By \cite{Gall-1} it is follows that $\delta_{i}=0$  for $i>30.$

 Summarizing the above results we get
\begin{te}
The  system of the following 147 covariants
$$
\begin{array}{l}
\kappa^{-1}(t)\\
\kappa^{-1}(dv_{i}),i=1,\ldots,3 , \kappa^{-1}(tr_{i}),i=1,\ldots,6, \kappa^{-1}(ch_{i}),i=1,\ldots,8, , \kappa^{-1}(pt_{i}),i=1,\ldots,10\\
 \kappa^{-1}(sh_{i}),i=1,\ldots,10, \kappa^{-1}(si_{i}),i=1,\ldots,12 , \kappa^{-1}(vi_{i}),i=1,\ldots,13 , \kappa^{-1}(de_{i}),i=1,\ldots,11\\
 \kappa^{-1}(des_{i}),i=1,\ldots,9 , \kappa^{-1}(odn_{i}),i=1,\ldots,8 , \kappa^{-1}(dvn_{i}),i=1,\ldots,13 , \kappa^{-1}(tryn_{i}),i=1,\ldots,9\\
 \kappa^{-1}(chot_{i}),i=1,\ldots,6 , \kappa^{-1}(ptn_{i}),i=1,\ldots,4 , \kappa^{-1}(shis_{i}),i=1,\ldots,5 , \kappa^{-1}(simn_{i}),i=1,\ldots,2 \\
vis_{i},i=1,\ldots,9,  \kappa^{-1}(devn), dvad,dvdv_1,dvdv_2, \kappa^{-1}(dvtr), dvsh, trd,
\end{array}
$$
is a complete system of the covariants for the  binary form of  degree 7.

\end{te}

\section{Appendix}

The degree-order distribution of $\overline C_8.$
\newpage

\begin{table}
\[ \qquad \qquad \text{order} \] \[ \text{degree} \; \; \begin{tabular}{|c||c|c|c|c|c|c|c|c|c|c|c|c|c|c|c|c|c|c|c|c|c|c|}
 \hline       &  0\phantom{0} & 1\phantom{0} &  2\phantom{0}  & 3\phantom{0}  & 4\phantom{0}   &5\phantom{0}  & 6\phantom{0}  &  7\phantom{0}  & 8 \phantom{0}   & 9\phantom{0}  &10 &11&12&13&14&15   \\ 
\hline 
\hline
 1                   & {} & {} & {} & {}  &  &  {} & 1    &  {} & {}   & {}  &&&&&&      \\ 
\hline
 2                   &  & {} & 1  &   &   &   &  1  &   & {}  & {}  &1&&&&&      \\ 

 \hline 3         & & {} &  &1  &   & 1  &    &  1 & {}   &1  &&1&&&&1     \\ 

\hline
 4                   &1 & {} &  & &  2 &    & 1   &  & 2  &   &1&&&&1&  \\ 
   
\hline
5                    &  & 1 &   & 2 &   & 2 & {}    & 2  & {}   & 2 &&&&1&&      \\  

 \hline
6                    &   & & 3  &   & 2  &  &2   &  {} & 2   & {}   &&&1&&&    \\   

\hline
7                    &  &3  &  &2    & {}  &  4 & {}    & 2 & {}   & {}  &&1&&&&   \\ 
 
 \hline
8                    & 3 &   & 3& {}  &  3  &  {}    &  3 & {}   & {}   &  &1&&&&&  \\    
                                                                            
 \hline
 9                   & &3  &  & 5  & {}  &  2 & {}    &  {} & {}   & 1  &&&&&&    \\ 

\hline
 10                   & &  & 4 & {}  & 4 &  {} & {}    &  {} & 1   & {}   &&&&&&    \\ 

\hline
 11                   & {} &5  & {} & 3  & {}  &  {} & {}    &  1 & {}   & {}  &&&&&&   \\ 

\hline
 12                   &6&  & 6 &  & {}  &  {} & 1   &  {} & {}   &    &&&&&&  \\ 

\hline
 13                   & {} &7 & {} & 1  & {}  &  1 & {}    &  {} & {}   &    &&&&&&  \\ 

\hline
 14                   &4 & {} &  & {}  & 2  &  {} & {}    &  {} & {}   & {}  &&&&&&      \\ 
\hline
 15                   &  & 3 &   &1 &   &  {} &    &  {} & {}  & {}  &&&&&&      \\ 

 \hline 16        & 2 &&3  &  &   &   &    &   & {}   &  &&&&&&     \\ 

\hline
17                  & & 2 &  & &   &    &    &  & {}   &   &&&&&&  \\ 
   
\hline
18                    & 9 &  &   &  &   &  & {}    &   & {}   & {} &&&&&&      \\  

 \hline
19                    &   &1 &   &   &   &  & {}    &  {} & {}   & {}   &&&&&&    \\   

\hline
20                   & 1 &  &  &    & {}  &  {} & {}    &  {} & {}   & {}  &&&&&&   \\ 
 
 \hline
21                   &  &   & & {}  &  {}  &  {}    &  {} & {}   & {}   &  &&&&&&  \\    
                                                                            
 \hline
 22                  & 2&  &  & {}  & {}  &  {} & {}    &  {} & {}   & {}  &&&&&&    \\ 

\hline
 23                   & &1  & {} & {}  & {}  &  {} & {}    &  {} & {}   & {}   &&&&&&    \\ 

\hline
24                   & {} &  & {} & {}  & {}  &  {} & {}    &  {} & {}   & {}  &&&&&&   \\ 

\hline
25                 & {} & 1  & {} & {}  & {}  &  {} & {}    &  {} & {}   &    &&&&&&  \\ 

\hline
26                  & 1 &{} & {} & {}  & {}  &  {} & {}    &  {} & {}   &    &&&&&&  \\

\hline
30                 & 1 &  & {} & {}  & {}  &  {} & {}    &  {} & {}   &    &&&&&&  \\ 

\hline \end{tabular} \] 
\end{table}

\end{document}